\documentclass[letterpaper, 10 pt, conference]{ieeeconf}
\IEEEoverridecommandlockouts
\usepackage[space]{cite}

\usepackage[tbtags]{amsmath}
\usepackage{amsmath}
\usepackage{amsfonts}
\usepackage{amssymb}
\usepackage{amsthm}
\usepackage{mathrsfs}
\usepackage{xspace}
\usepackage{xparse}
\usepackage{mathtools}

\let\labelindent\relax
\usepackage[inline]{enumitem}

\usepackage{algorithm}
\usepackage{algpseudocode}
\usepackage{setspace}

\usepackage{tabularx}
\usepackage{booktabs}
\usepackage{multirow}

\usepackage{tikz}
\usetikzlibrary{positioning}

\usepackage{epstopdf}
\usepackage{color, soul}

\usepackage{hyperref}

\newtheorem{defn}{Definition}

\newtheorem{rem}{Remark}

\newtheorem{assum}{Assumption}

\begin{document}

\title{
Stochastic Optimal Control via Hilbert Space Embeddings of Distributions
}

\author{
Adam~J.~Thorpe,~\IEEEmembership{Student~Member,~IEEE,}
Meeko~M.~K.~Oishi,~\IEEEmembership{Senior Member,~IEEE}%
\thanks{%
    This material is based upon work supported by the National Science Foundation under NSF Grant Number CNS-1836900.  Any opinions, findings, and conclusions or recommendations expressed in this material are those of the authors and do not necessarily reflect the views of the National Science Foundation.
    The NASA University Leadership initiative (Grant \#80NSSC20M0163) provided funds to assist the authors with their research, but this article solely reflects the opinions and conclusions of its authors and not any NASA entity.
    This research was supported in part by the Laboratory Directed Research and Development program at Sandia National Laboratories, a multimission laboratory managed and operated by National Technology and Engineering Solutions of Sandia, LLC., a wholly owned subsidiary of Honeywell International, Inc., for the U.S. Department of Energy’s National Nuclear Security Administration under contract DE-NA-0003525.  The views expressed in this article do not necessarily represent the views of the U.S. Department of Energy or the United States Government.
}
\thanks{A. Thorpe and M. Oishi are with Electrical
  \& Computer Eng., University of New Mexico, Abq., NM. 
  Email: {\tt\{ajthor,oishi\}@unm.edu}.
}
}

\maketitle

\begin{abstract}
Kernel embeddings of distributions have recently gained significant attention in the machine learning community as a data-driven technique for representing 
probability distributions. Broadly, these techniques enable efficient computation of expectations by representing integral operators as elements in a reproducing kernel Hilbert space. We apply these techniques to the area of stochastic optimal control theory and present a method to compute approximately optimal policies for stochastic systems with arbitrary disturbances. Our approach reduces the optimization problem to a linear program, which can easily be solved via the Lagrangian dual, without resorting to gradient-based optimization algorithms. We focus on discrete-time dynamic programming, and demonstrate our proposed approach on a linear regulation problem, and on a nonlinear target tracking problem. This approach is broadly applicable to a wide variety of optimal control problems, and provides a means of working with stochastic systems in a data-driven setting.
\end{abstract}


\section{Introduction}


Stochastic systems are ubiquitous, however 
most methods for control of stochastic systems are reliant upon accurate modeling not only of the dynamics, but also of the stochastic processes of the system.
As autonomous systems become commonplace, and direct human interaction with autonomy become more pervasive, presumptions of linearity and Gaussian stochasticity become questionable, as they could lead to control solutions that are confusing, non-intuitive, or simply incorrect. Robust solutions that work well when uncertainty is bounded by known values may be excessively conservative, and cannot accommodate long-tail phenomena.  
In contrast, data-driven approaches do not rely upon any prior assumptions on the dynamics or stochasticity of the system. Data-driven approaches 
have garnered considerable interest recently, due to the capabilities of learning algorithms to handle systems with nonlinear dynamics and unknown disturbances.

We propose a method 
for data-driven controller synthesis 
based on \emph{conditional distribution embeddings} \cite{song2009hilbert}, a nonparametric learning technique that uses a sample of system observations to construct a model of the stochastic system dynamics as an element in a high-dimensional Hilbert space of functions known as a \emph{reproducing kernel Hilbert space}. 
These techniques leverage functional analysis and statistical learning theory to empirically estimate the stochastic kernel using data. 
As a nonparametric technique, kernel methods are inherently data-driven, and do not rely upon prior assumptions placed upon the data or exploit system structure.
These techniques have been applied to Markov models \cite{grunewalder2012modelling}, partially-observable systems \cite{nishiyama2012hilbert, song2010hilbert}, and more recently, to robust optimization approaches \cite{zhu2020kernel, zhu2020kernel_b}.
Furthermore, these techniques admit finite sample bounds which show convergence in probability as the number of samples tend to infinity \cite{song2009hilbert}. 
%
A Hilbert space framework is particularly well-suited to stochastic optimal control problems \cite{luenberger1997optimization}, primarily because Hilbert spaces are a generalization of inner product spaces to an infinite-dimensional setting, meaning they encompass many optimization problems of interest (note that $\mathbb{R}^{n}$ is a Hilbert space).


The use of kernel methods for policy synthesis is well-motivated in literature, especially in the area of reinforcement learning (RL) \cite{ormoneit2002kernel}.
Methods have been developed to optimize a policy in an RKHS via functional gradient descent \cite{bagnell2003policy}, \cite{lever2015modelling}.  
Other approaches rely upon value iteration, approximating the value function as an intermediate step in order to compute an optimal control input \cite{grunewalder2012modelling}.  
However, most of these approaches face significant computational challenges due to the sampling schemes used by RL, the need for knowledge of a gradient, or reliance upon iterative numerical methods. 
Some progress has been made to alleviate these issues, for example using stochastic factorization \cite{barreto2016practical}.

Our main contribution is a data-driven algorithm for computing approximately optimal policies for arbitrary discrete-time stochastic dynamical systems. 
The novelty of our approach is the use of conditional distribution embeddings to formulate an optimal control problem as a linear program within a reproducing kernel Hilbert space, which can be solved efficiently via the Lagrangian dual. 
Our approach is model-free, since it relies only upon data collected from prior observations of the system execution, meaning that 
it is amenable to systems with arbitrary disturbances and nonlinear dynamics. 
Because we formulate the optimal control problem as a linear program, we do not rely upon gradient-based algorithms to compute an optimal solution, and thus do not impose a specific structure on the policy for the purpose of computing a functional gradient. 
The main difficulty associated with our approach is the dependence of the computational complexity on sample size (generally $\mathcal{O}(M^{3})$), as with all kernel based approaches.
This arises from the presence of a matrix inverse operation and the large number of observations needed to fully characterize the stochasticity of a system. 
Fortunately, numerous approaches to reducing the computational burden of kernel methods have been explored, such as \cite{rahimi2007random, le2013fastfood}, which use Fourier transforms and Gaussian matrix approximations to reduce the computational complexity to log-linear time.

The paper is structured as follows:
In section \ref{section: preliminaries}, we define the problem and describe the preliminary theory of embedding distributions in reproducing kernel Hilbert spaces in section \ref{section: embedding stochastic kernels}. 
We then present our method in section \ref{section: hilbert space policy optimization} to compute the optimal policy and present an extension of our proposed approach to solve dynamic programming problems over a finite time horizon.
In section \ref{section: numerical results}, we demonstrate our proposed approach on a simple stochastic integrator system for the purpose of validation against a known result, and then on a target tracking problem using nonlinear, nonholonomic vehicle dynamics.
Concluding remarks are presented in section \ref{section: concluding remarks}.


\section{Preliminaries}
\label{section: preliminaries}

We use the following notation throughout:
Let $E$ be an arbitrary nonempty space, and denote the $\sigma$-algebra on $E$ by $\mathcal{E}$.
%
If $E$ is a topological space \cite{ccinlar2011probability},
the $\sigma$-algebra generated by the set of all open subsets of $E$
is called the Borel $\sigma$-algebra, denoted by $\mathscr{B}(E)$.
%
Let $(\Omega, \mathcal{F}, \mathbb{P})$ denote a probability space,
where $\mathcal{F}$ is the $\sigma$-algebra on $\Omega$ and
$\mathbb{P} : \mathcal{F} \rightarrow [0, 1]$ is a \emph{probability measure}
on the measurable space $(\Omega, \mathcal{F})$.
A measurable function $X : \Omega \rightarrow E$ is called a \emph{random variable} taking values in $(E, \mathcal{E})$. The image of $\mathbb{P}$ under $X$, $\mathbb{P}(X^{-1}A)$, $A \in \mathcal{E}$ is called the \emph{distribution} of $X$.
Let $T$ be an arbitrary set, and for each $t \in \mathcal{T}$, let $X_{t}$ be a random variable.
The collection of random variables $\lbrace X_{t} : t \in \mathcal{T} \rbrace$
on $(\Omega, \mathcal{F})$ is a \emph{stochastic process}.
We define a \emph{stochastic kernel} according to \cite{ccinlar2011probability}.

\begin{defn}[Stochastic Kernel]
    Let $(E, \mathcal{E})$ and $(F, \mathcal{F})$ be measurable spaces with $\sigma$-algebras $\mathcal{E}$ and $\mathcal{F}$, respectively.
    A stochastic kernel is a map $\kappa : \mathcal{F} \times E \to [0, 1]$, where:
    \begin{enumerate*}[mode=unboxed]
        \item 
        $x \mapsto \kappa(B \,|\, x)$ is $\mathcal{E}$-measurable for all $B \in \mathcal{F}$;
        
        \item 
        $B \mapsto \kappa(B \,|\, x)$ is a probability measure on $(F, \mathcal{F})$ for all $x \in E$.
    \end{enumerate*}
\end{defn}


\subsection{System Model}

Consider a Markov control process, which is defined in \cite{puterman2005markov} as a $3$-tuple, $(\mathcal{X}, \mathcal{U}, Q)$, consisting of:
\begin{itemize}
    \item 
    A Borel space $\mathcal{X} \subseteq \mathbb{R}^{n}$ called the state space;
    \item 
    A compact Borel space $\mathcal{U} \subset \mathbb{R}^{m}$ called the control space; and
    \item 
    A stochastic kernel $Q : \mathscr{B}(\mathcal{X}) \times \mathcal{X} \times \mathcal{U} \to [0, 1]$ that assigns a probability measure $Q(\cdot \,|\, x, u)$ to each $(x, u) \in \mathcal{X} \times \mathcal{U}$ on the measurable space $(\mathcal{X}, \mathscr{B}(\mathcal{X}))$.
\end{itemize}

The system evolves from an initial condition $x_{0} \in \mathcal{X}$, which may be chosen from an initial distribution $\mathbb{P}_{0}$ on $\mathcal{X}$, over a finite time horizon $t = 0, 1, \ldots, N$, $N \in \mathbb{N}_{+}$.
As the system evolves, the control actions $u_{0}, u_{1}, \ldots, u_{N-1}$ are chosen from a Markov control policy $\pi$.

\begin{defn}[Markov Policy, {\cite[Definition 8.2]{bertsekas1978stochastic}}]
    \label{def: markov policy}
    A Markov policy $\pi$ is a sequence $\pi = \lbrace \pi_{0}, \pi_{1}, \ldots \pi_{N-1} \rbrace$ of universally measurable stochastic kernels, where for each $t = 0, 1, \ldots, N-1$, the stochastic kernel $\pi_{t} : \mathscr{B}(\mathcal{U}) \times \mathcal{X} \to [0, 1]$ assigns a probability measure $\pi_{t}(\cdot \,|\, x)$ to every $x \in \mathcal{X}$ on the measurable space $(\mathcal{U}, \mathscr{B}(\mathcal{U}))$.
\end{defn}


\subsection{Problem Formulation}

We assume that the stochastic kernel $Q$ is unknown, but that a sample of observations of the system evolution is available. 

\begin{assum}
    \label{assum: stochastic kernel unknown}
    We assume that $Q$ is unknown, but that a sample $\mathcal{S} = \lbrace (x_{i}, u_{i}, x_{i}{}^{\prime}) \rbrace_{i=1}^{M}$ of size $M \in \mathbb{N}_{+}$ is available, where $x_{i}{}^{\prime} \sim Q(\cdot \,|\, x_{i}, u_{i})$ and $u_{i}$ is selected randomly from the set of admissible control inputs.
\end{assum}

Consider an arbitrary cost function $c : \mathcal{X} \to \mathbb{R}$, which we assume is a continuous, bounded functional that lies in a Hilbert space of functions $\mathscr{H}$.
At any time instant $t$, we seek to minimize $c$ by selecting the \emph{distribution} $\pi_{t}$ on $(\mathcal{U}, \mathscr{B}(\mathcal{U}))$ which minimizes the following unconstrained minimization problem:
\begin{align}
    \label{eqn: optimization problem}
    \begin{aligned}[c]
        \min_{\pi_{t}} \quad
        & J_{t}(\pi) = \int_{\mathcal{U}} \int_{\mathcal{X}} c(y) Q(\mathrm{d} y \,|\, x, v) \pi_{t}(\mathrm{d} v \,|\, x)
    \end{aligned}
\end{align}

The primary difficulty in solving \eqref{eqn: optimization problem} is 
that without knowledge of $Q$, the integral in \eqref{eqn: optimization problem} is intractable.
Thus, we seek to form an approximate optimization problem by approximating the integral in \eqref{eqn: optimization problem} using a sample $\mathcal{S}$ taken i.i.d. from $Q$ as an element in a Hilbert space of functions.  
By optimizing the approximate problem, we obtain an approximate solution. 
Thus, we additionally seek to ensure that the approximate optimization problem converges in probability to the true optimization problem as the sample size increases.

According to \cite{bertsekas1978stochastic}, in most cases, the optimal Markov policy for a system can be viewed as nonrandomized, or deterministic, meaning the stochastic kernel assigns a probability measure with mass one at a single element in $\mathcal{U}$ to each $x \in \mathcal{X}$. 
According to \cite{bertsekas1978stochastic, puterman2005markov}, the set of nonrandomized policies is a subset of the set of all randomized policies, meaning we can search among the class of randomized policies in Hilbert space to find an optimal policy which minimizes \eqref{eqn: optimization problem}. 


\section{Embedding Stochastic Kernels in an RKHS}
\label{section: embedding stochastic kernels}

Let $\mathscr{H}$ be a Hilbert space of functions of the form $\mathcal{X} \to \mathbb{R}$ with inner product $\langle \cdot, \cdot \rangle_{\mathscr{H}}$ and the induced norm $\lVert \cdot \rVert_{\mathscr{H}}$.

\begin{defn}[RKHS, {\cite{aronszajn1950theory}}]
     A Hilbert space $\mathscr{H}$ is a reproducing kernel Hilbert space (RKHS) if there exists a positive definite \cite[Definition~4.12]{steinwart2008support} kernel function $k$ that satisfies the following properties:
     \begin{align}
         & k(x, \cdot) \in \mathscr{H}, && \forall x \in \mathcal{X} \\
         \label{eqn: reproducing property}
         & f(x) = \langle f, k(x, \cdot) \rangle_{\mathscr{H}}, && \forall f \in \mathscr{H}, x \in \mathcal{X}
     \end{align}
     where \eqref{eqn: reproducing property} is known as the reproducing property, and for any $x, x' \in \mathcal{X}$, we denote $k(x, \cdot) \in \mathscr{H}$ as a function on $\mathcal{X}$ such that $x' \mapsto k(x, x')$.
\end{defn}

    

\begin{rem}
    Alternatively, by the Moore-Aronszajn theorem \cite{aronszajn1950theory}, we can define an RKHS by first specifying a kernel $k$ and obtain a corresponding RKHS as the closure of the span of kernel functions. 
\end{rem}

Given $(x, u) \in \mathcal{X} \times \mathcal{U}$, let $Q(\cdot \,|\, x, u)$ be a conditional probability measure on $\mathcal{X}$.
According to \cite{song2009hilbert}, if the following sufficient condition holds:
\begin{equation}
    \label{eqn: sufficient condition}
    \int_{\mathcal{X}} \sqrt{k(y, y)} Q(\mathrm{d} y \,|\, x, u) < \infty
\end{equation}
then there exists an element $m(x, u) \in \mathscr{H}$ called a \emph{conditional distribution embedding}, where
\begin{equation}
    m(x, u) := \int_{\mathcal{X}} k(y, \cdot) Q(\mathrm{d} y \,|\, x, u)
\end{equation}
By the reproducing property of $k$ in $\mathscr{H}$, for any $f \in \mathscr{H}$, we can evaluate the integral with respect to $Q(\cdot \,|\, x, u)$ as an inner product with the embedding $m(x, u)$:
\begin{align}
    \label{eqn: embedding reproducing property}
    \langle f, m(x, u) \rangle_{\mathscr{H}} 
    &= \biggl\langle f, \int_{\mathcal{X}} k(y, \cdot) Q(\mathrm{d} y \,|\, x, u) \biggr\rangle_{\mathscr{H}} \\
    &= \int_{\mathcal{X}} \langle f, k(y, \cdot) \rangle_{\mathscr{H}} Q(\mathrm{d} y \,|\, x, u) \\
    &= \int_{\mathcal{X}} f(y) Q(\mathrm{d} y \,|\, x, u)
\end{align}
Intuitively, the element $m \in \mathscr{H}$ corresponds to the dynamics of the system at the point $(x, u)$.
In other words, if the integral exists, then we can embed the integral operator with respect to the probability measure in $\mathscr{H}$ and evaluate the integral via the reproducing property of $k$ in $\mathscr{H}$. 

However, in a data-driven setting, the stochastic kernel $Q$ is unknown, which means the embedding $m(x, u)$ is also unknown. Instead, we can empirically estimate the stochastic kernel using a sample of observations taken from $Q$. 


\subsection{Empirical Embeddings Using Observations}




Consider a sample $\mathcal{S} = \lbrace (x_{i}, u_{i}, x_{i}{}^{\prime}) \rbrace_{i=1}^{M}$ of size $M \in \mathbb{N}_{+}$, taken i.i.d. from $Q$, where $x_{i}$ and $u_{i}$ are taken randomly from the state and control spaces $\mathcal{X}$ and $\mathcal{U}$, respectively, and $x_{i}{}^{\prime} \sim Q(\cdot \,|\, x, u)$. 
As shown in \cite{grunewalder2012conditional}, we can compute an empirical estimate $\hat{m}$ of $m$ as the solution to a regularized least-squares problem, given by:
\begin{equation}
    \label{eqn: regularized least squares problem}
    \min_{\hat{m}} \frac{1}{M} \sum_{i=1}^{M} \lVert k(x_{i}{}^{\prime}, \cdot) - \hat{m}(x_{i}, u_{i}) \rVert_{\mathscr{H}}^{2} + \lambda \lVert \hat{m} \rVert_{\mathscr{Q}}^{2}
\end{equation}
where $\lambda > 0$ is the regularization parameter and $\mathscr{Q}$ is a \emph{vector-valued} RKHS \cite{grunewalder2012conditional}. 
As shown in \cite{micchelli2005learning, grunewalder2012conditional}, by the representer theorem, the solution $\hat{m}$ to \eqref{eqn: regularized least squares problem} is unique and has the following form:
\begin{equation}
    \label{eqn: estimate form}
    \hat{m}(x, u) = \sum_{i=1}^{M} \beta_{i}(x, u) k(x_{i}{}^{\prime}, \cdot)
\end{equation}
where $\beta(x, u) \in \mathbb{R}^{M}$ is a vector of real-valued coefficients that depends on the conditioning variables $x$ and $u$. 
The problem in \eqref{eqn: regularized least squares problem} admits a closed-form solution, given by:
\begin{equation}
    \label{eqn: estimate closed form solution}
    \hat{m}(x, u) = \Phi^{\top} (\Psi \Psi^{\top} + \lambda M I)^{-1} \Psi k(x, \cdot) k(u, \cdot)
\end{equation}
where $\Phi$ and $\Psi$ are called \emph{feature vectors}, with elements given by $\Phi_{i} = k(x_{i}{}^{\prime}, \cdot)$ and $\Psi_{i} = k(x_{i}, \cdot) k(u_{i}, \cdot)$, respectively.
For simplicity, we denote $W = (\Psi \Psi^{\top} + \lambda M I)^{-1}$ and let $\beta(x, u) = W \Psi k(x, \cdot) k(u, \cdot)$, such that $\hat{m}(x, u) = \Phi^{\top} \beta(x, u)$. 
Using $\hat{m}(x, u)$, we can approximate the expectation with respect to $Q(\cdot \,|\, x, u)$ for any $f \in \mathscr{H}$ as:
\begin{equation}
    \langle f, \hat{m}(x, u) \rangle_{\mathscr{H}} \approx \int_{\mathcal{X}} f(y) Q(\mathrm{d} y \,|\, x, u)
\end{equation}

Further, if the kernel function $k$ is \emph{universal} \cite{micchelli2006universal}, then the embedding is injective, meaning there exists a unique representation of the distribution in $\mathscr{H}$. 
In short, a universal kernel $k$ allows us to approximate any arbitrary real-valued function using \eqref{eqn: estimate form} arbitrarily well as the number of samples tends to infinity.
A commonly used kernel function which satisfies this property is the Gaussian kernel $k(x, x') = \exp(-\lVert x - x' \rVert_{2}^{2}/2\sigma^2)$, $\sigma > 0$.
Additionally, the estimate $\hat{m}$ converges in probability to the true embedding $m$ as the number of samples $M$ tends to infinity and $\lambda \to 0$ \cite{song2009hilbert, song2010nonparametric}.
This means the estimate $\hat{m}$ is a \emph{consistent} estimator of the true embedding, and the integral of a function $f \in \mathscr{H}$ with respect to $Q$ converges in probability to the true result as the sample size increases. 


\section{Policy Optimization in Hilbert Space}
\label{section: hilbert space policy optimization}

Consider the problem in \eqref{eqn: optimization problem} where $c \in \mathscr{H}$.
As shown in \cite{song2009hilbert}, if the sufficient condition in \eqref{eqn: sufficient condition} holds, then there exists a conditional distribution embedding $m(x, u)$ such that for any $c \in \mathscr{H}$,
\begin{equation}
    \label{eqn: expected cost}
    \int_{\mathcal{X}} c(y) Q(\mathrm{d} y \,|\, x, u)
    = \langle c, m(x, u) \rangle_{\mathscr{H}}
\end{equation}
This allows us to evaluate the expected cost at a particular $(x, u) \in \mathcal{X} \times \mathcal{U}$ as an inner product in Hilbert space.
Let $\pi$ be a Markov policy as in Definition \ref{def: markov policy}.
Taking the integral of \eqref{eqn: expected cost} with respect to the Markov policy $\pi$, we obtain the objective function $J_{t}(\pi)$ in \eqref{eqn: optimization problem}. 
By linearity of the integral and the inner product, we can rewrite the objective using the inner product in \eqref{eqn: expected cost} to obtain:
\begin{align}
    \label{eqn: expected cost under polict pi}
        J_{t}(\pi) 
        &= \biggl\langle c, \int_{\mathcal{U}} m(x, v) \pi_{t}(\mathrm{d} v \,|\, x) \biggr\rangle_{\mathscr{H}}
\end{align}
where the integral term on the right hand side of \eqref{eqn: expected cost under polict pi} can be interpreted as a representation in $\mathscr{H}$ of the closed-loop dynamics under a policy $\pi$.

However, the integral in \eqref{eqn: expected cost under polict pi} is intractable, since according to Assumption \ref{assum: stochastic kernel unknown}, the stochastic kernel $Q$ (and thus the embedding $m$) is unknown. 
Instead, we compute an empirical estimate $\hat{m}$ of $m$ using a sample $\mathcal{S}$ taken i.i.d. from $Q$.
Recall from \eqref{eqn: estimate closed form solution} that the empirical estimate has the form $\hat{m}(x, u) = \Phi^{\top} W \Psi k(x, \cdot) k(u, \cdot)$.
We then substitute the estimate for the true embedding to approximate the integral in \eqref{eqn: expected cost under polict pi}.
\begin{align}
    \label{eqn: approximate expected cost under policy pi}
    \begin{split}
        \int_{\mathcal{U}} m(x, v) \pi_{t}(\mathrm{d} v \,|\, x) 
        &\approx \int_{\mathcal{U}} \hat{m}(x, v) \pi_{t}(\mathrm{d} v \,|\, x) \\
        &= \Phi^{\top} W \Psi k(x, \cdot) \int_{\mathcal{U}} k(v, \cdot) \pi_{t}(\mathrm{d} v \,|\, x)
    \end{split}
\end{align}
Recall that the policy is a collection of stochastic kernels indexed by time, which means that at a given time $t$, the policy can be represented by a conditional distribution embedding. 
Thus, it is natural to consider the policy at a particular time as a collection of elements in an RKHS parameterized by $x \in \mathcal{X}$, which admits a representation in terms of finite support.
Let $\lbrace \tilde{u}_{j} \rbrace_{j=1}^{P}$ be a collection of admissible control inputs.
We propose the following representation for the policy $\pi_{t}$:
\begin{equation}
    \label{eqn: policy representation}
    \hat{p}_{t}(x) = \sum_{j=1}^{P} \alpha_{j}(x) k(\tilde{u}_{j}, \cdot)
\end{equation}
where $\alpha(x) \in \mathbb{R}^{P}$ is a vector of real-valued coefficients that depends on $x \in \mathcal{X}$. 
Using \eqref{eqn: policy representation}, we can approximate \eqref{eqn: approximate expected cost under policy pi} as:
\begin{equation}
    \Phi^{\top} W \Psi k(x, \cdot) \int_{\mathcal{U}} k(v, \cdot) \pi_{t}(\mathrm{d} v \,|\, x)
    \approx \Phi^{\top} W \Psi k(x, \cdot) \Upsilon^{\top} \alpha(x)
\end{equation}
where $\Upsilon$ is a feature vector with elements $\Upsilon_{j} = k(\tilde{u}_{j}, \cdot)$.
Thus, we can approximate the objective function $J_{t}(\pi)$ by:
\begin{equation}
    \label{eqn: unbounded optimization problem}
    \int_{\mathcal{U}} \int_{\mathcal{X}} c(y) Q(\mathrm{d} y \,|\, x, v) \pi_{t}(\mathrm{d} v \,|\, x)
    \approx \boldsymbol{c}^{\top} W \Psi k(x, \cdot) \Upsilon^{\top} \alpha(x)
\end{equation}
where $\boldsymbol{c}$ is a vector with elements $\boldsymbol{c}_{i} = c(x_{i}{}^{\prime})$.
Thus, we form an approximation of the objective in \eqref{eqn: optimization problem}, which converges in probability to the true optimization problem as the sample sizes $M$ and $P$ increase \cite{song2009hilbert}. 
Now, instead of minimizing over the distribution $\pi_{t}$, we can view the approximate optimization problem as finding $\alpha(x) \in \mathbb{R}^{P}$ which minimizes \eqref{eqn: unbounded optimization problem}. 
However, minimizing $\alpha(x)$ in \eqref{eqn: unbounded optimization problem} is unbounded below, which makes the problem unsolvable.
As such, additional constraints are required to ensure that the problem admits a feasible solution.
Note that intuitively, $\hat{p}_{t}(x)$ is an approximation of the distribution $\pi_{t}(\cdot \,|\, x)$ at time $t$. 
Because of this, we can view the coefficients $\alpha(x)$ as a vector of probabilities which weight the nonlinear transformations of the control inputs. 
Thus, we place additional constraints on the coefficients $\alpha(x)$, and form the approximate optimization problem, constraining the values of $\alpha(x)$ such that they are non-negative and sum to one:
\begin{subequations}
    \label{eqn: approximate optimization problem}
    \begin{align}
        \min_{\alpha(x) \in \mathbb{R}^{P}} \quad 
        & \boldsymbol{c}^{\top} W \Psi k(x, \cdot) \Upsilon^{\top} \alpha(x) \\
        \textnormal{s.t.} \quad 
        & \sum_{j=1}^{P} \alpha_{j}(x) = 1 \\
        & 0 \preceq \alpha(x)
    \end{align}
\end{subequations}
Since $\alpha(x)$ is an indirect weighting on the inputs that depends on the state $x$, we interpret $\alpha(x)$ as a probability weighting of the control inputs $\tilde{u}_{j}$.
Note that \eqref{eqn: approximate optimization problem} is a linear program in standard form \cite{boyd2004convex}, and that we can solve \eqref{eqn: approximate optimization problem} via the Lagrangian dual. Let $\nu \in \mathbb{R}$ be a dual variable, and for simplicity, let $C(x) = \boldsymbol{c}^{\top} W \Psi k(x, \cdot) \Upsilon^{\top}$. 
The dual problem is given by:
\begin{subequations}
    \label{eqn: dual problem}
    \begin{align}
        \max \quad 
        & -\nu \\
        \textnormal{s.t.} \quad 
        & -\boldsymbol{1} \nu \preceq C(x)^{\top}
    \end{align}
\end{subequations}
where $\boldsymbol{1}$ is a vector of all ones.
From \cite[\S 4]{boyd2004convex}, \eqref{eqn: dual problem} has an optimal solution, given by $\min_{i} \lbrace C_{i}(x)^{\top} \rbrace$, which means the optimal solution $\alpha(x)^{*}$ to \eqref{eqn: approximate optimization problem} is a vector of all zeros, except $\alpha_{i}(x)^{*} = 1$.
In other words, we choose the control input that corresponds to the minimal value of $C(x)$. 


\subsection{Application to Approximate Dynamic Programming}
\label{section: approximate dynamic programming}

Many optimal control problems can be formulated as dynamic programs. 
Consider the following problem with an additive cost, in which we seek a policy $\pi$ that minimizes the following optimization problem \cite{bertsekas1978stochastic}:
    \begin{align}
        \label{eqn: standard dynamic programming problem}
        \min_{\pi} \quad 
        & J_{N}(\pi) = \mathbb{E}_{\pi} \biggl[ g_{N}(x_{N}) + \sum_{i=0}^{N-1} g_{t}(x_{t}, u_{t}) \biggr]
    \end{align}
where $N \in \mathbb{N}_{+}$ is the time horizon, $\pi$ is the control policy, $g_{N}$ is the terminal cost for ending in state $x_{N}$, $g_{t}$ is the cost at time $t$ of taking action $u_{t} \sim \pi_{t}(\cdot \,|\, x_{t})$ while in state $x_{t}$, and the expectation is uniquely determined by the initial distribution $\mathbb{P}_{0}$ and the Markov policy $\pi$ (see \cite[Definition 8.3]{bertsekas1978stochastic} for more details). 
The problem in \eqref{eqn: standard dynamic programming problem} can be rewritten via the Chapman-Kolmogorov identity and the Markov property as a sequence of sub-problems, where the problem is solved backward in time via \emph{backward recursion} \cite{bertsekas1978stochastic}.
We define the \emph{value functions} $V_{t} : \mathcal{X} \to \mathbb{R}$ for all $t = 0, 1, \ldots, N-1$ as:
\begin{equation}
    \label{eqn: value function recursion}
    V_{t}(x) = \max_{\pi_{t}} \int_{\mathcal{U}} \int_{\mathcal{X}} g(x, v) + V_{t+1}(y) Q(\mathrm{d} y \,|\, x, v) \pi_{t}(\mathrm{d} v \,|\, x)
\end{equation}
initialized with $V_{N}(x_{N}) = g_{N}(x_{N})$.
Then the solution to \eqref{eqn: standard dynamic programming problem} is equivalent to solving a sequence of sub-problems given by \eqref{eqn: value function recursion}, and iteratively substituting the solutions into the subsequent value function.

We can apply \eqref{eqn: approximate optimization problem} in this context to solve for the optimal control policy when the dynamics and stochasticity are not known, but a sample $\mathcal{S}$ is available. 
In this case, we solve \eqref{eqn: value function recursion} at each time step $t$ using \eqref{eqn: approximate optimization problem}.
By approximating and recursively substituting the solution to \eqref{eqn: value function recursion} into the subsequent value function, we obtain an approximately optimal control policy $\pi^{*} \approx \arg \min_{\pi} J_{N}(\pi)$ which approximately minimizes the cost in \eqref{eqn: standard dynamic programming problem}. 

This means we can compute the approximately optimal policy for a problem without exploiting knowledge of the system dynamics or the structure of the disturbance. By solving for the approximately optimal policy using \eqref{eqn: approximate optimization problem}, we avoid intractable integrals in the stochastic optimal control problem and can compute the policy as a linear operation in a Hilbert space of functions.
Additionally, this approach is largely agnostic to the dimensionality of the system, since the system dimensionality only directly affects the computation of the kernel function. For example, the Gaussian kernel scales linearly as the system dimensionality is increased. 
However, higher-dimensional systems typically require a larger sample size in order to fully characterize the dynamics of the system, which can be computationally prohibitive if the sample size is large. 






\section{Numerical Results}
\label{section: numerical results}

We demonstrate our approach on a $2$-D discrete-time stochastic integrator system for the purpose of verification, and on a target tracking problem with nonholonomic vehicle dynamics to demonstrate the utility of the approach. 
For all problems, we used a Gaussian kernel $k(x, x') = \exp(-\lVert x - x' \rVert_{2}^{2}/2 \sigma^{2})$.
Following \cite{song2009hilbert}, we chose the regularization parameter to be $\lambda = 1/M^{2}$, where $M$ is the sample size, as the default parameter for our calculations.
In practice, the parameters $\sigma$ and $\lambda$ are chosen via cross-validation, where $\sigma$ is selected according to the relative ``spacing'' of the observations, and $\lambda$ is a ``smoothness'' parameter chosen such that $\lambda \to 0$ as $M \to \infty$.
A more detailed discussion of parameter selection is outside the scope of the current work (see \cite{caponnetto2007optimal, song2009hilbert} for more information).
Numerical experiments were performed in Matlab on an AWS cloud computing instance, and computation times were obtained using Matlab’s Performance Testing Framework.

Code to reproduce the analysis and all figures is provided at:
\url{github.com/unm-hscl/ajthor-CDC2021}.


\subsection{Double Integrator System}

We consider the problem of regulation for a system whose dynamics are governed by a stochastic $2$-D discrete time stochastic integrator system, without any knowledge of the dynamics or the stochastic processes. 
That is, we seek a distribution $\pi$ which minimizes the following optimization problem:
\begin{align}
    \min_{\pi} \quad 
    & \int_{\mathcal{U}} \int_{\mathcal{X}} c(y) Q(\mathrm{d} y \,|\, x, v) \pi(\mathrm{d} v \,|\, x)
\end{align}
where $Q$ is a representation of the \emph{unknown} system dynamics as a stochastic kernel. 
For the purpose of comparison, we chose the cost function $c : \mathcal{X} \to \mathbb{R}$ to be the norm function:
\begin{equation}
    c(x) = \lVert x \rVert_{2}
\end{equation}
which serves to drive the system to the origin.
The dynamics for a $2$-D discrete time stochastic integrator system with sampling time $T_{s}$ are given by:
\begin{equation}
    \label{eqn: double integrator dynamics}
    x_{t+1} = 
    \begin{bmatrix}
        1 & T_{s} \\ 0 & 1
    \end{bmatrix}
    x_{t} + 
    \begin{bmatrix}
        T_{s}^{2}/2 \\ T_{s}
    \end{bmatrix}
    u_{t} + w_{t}
\end{equation}
where $x_{t} \in \mathcal{X}$ is the state, $u_{t} \in \mathcal{U}$ is the control input, which we specify to lie within the bounds $u_{t} \in [-1, 1]$, and $w$ is a stochastic process, comprised of the random variables $w_{t}$ on the measurable space $(\mathbb{R}^{p}, \mathscr{B}(\mathbb{R}^{p}))$. 
We consider three distributions for the disturbance:
\begin{enumerate*}[mode=unboxed]
    \item 
    A Gaussian distribution $w_{t} \sim \mathcal{N}(0, \Sigma)$, $\Sigma = 0.01 I$;
    
    \item 
    A beta distribution
    $w_{t} \sim 0.1 \mathrm{Beta}(\alpha, \beta)$, with a probability density function (PDF) given by:
    \begin{equation}
        f(x \,|\, \alpha, \beta) =
        \frac{\Gamma(\alpha + \beta)}{\Gamma(\alpha) \Gamma(\beta)}
        x^{\alpha-1} (1-x)^{\beta-1}
    \end{equation}
    where $\Gamma$ is the Gamma function and shape parameters $\alpha = 2$, $\beta = 0.5$; and
    
    \item
    An exponential distribution
    $w_{t} \sim 0.01 \mathrm{Exp}(\alpha)$,
    with $\alpha = 3$ and PDF $f(x \,|\, \alpha) = \alpha \exp (-\alpha x)$.
\end{enumerate*}

\begin{figure}[!t]
    \centering
    \includegraphics[width=240pt, height=240pt]{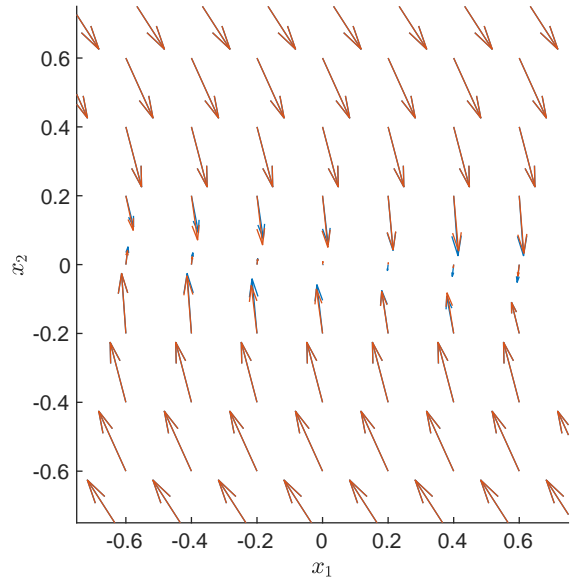}
    \caption{%
        Vector field showing the optimal closed-loop dynamics of a $2$-D discrete-time integrator system under an optimal control strategy computed via CVX (blue). The vector field of the closed-loop dynamics of a stochastic integrator system with a Gaussian disturbance computed using our proposed algorithm (orange).
    }
    \label{fig: double integrator}
\end{figure}

\begin{figure}[!t]
    \centering
    \includegraphics[width=120pt, height=120pt]{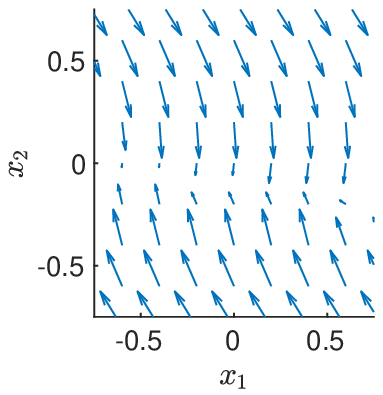} \hfill
    \includegraphics[width=120pt, height=120pt]{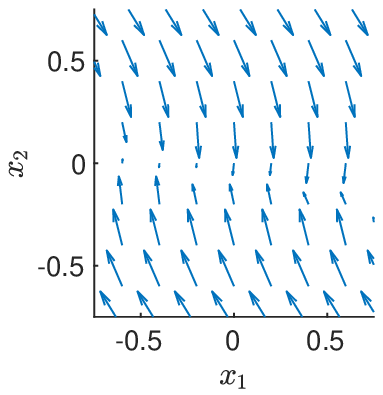}
    \caption{%
        (Left) Vector field of the approximately optimal closed-loop dynamics of a $2$-D stochastic integrator system with a Beta disturbance, where the approximately optimal policy is computed using our proposed approach.
        (Right) Vector field of the approximately optimal closed-loop system with an exponential disturbance. 
    }
    \label{fig: double integrator non gaussian}
\end{figure}

We consider a sample $\mathcal{S} = \lbrace (x_{i}, u_{i}, x_{i}{}^{\prime}) \rbrace_{i=1}^{M}$ of observations of size $M = 1600$ taken i.i.d. from $Q$, a representation of \eqref{eqn: double integrator dynamics} as a Markov control process. The states $x_{i} \in \mathcal{X}$ were selected uniformly in the range $x_{i} \in [-1, 1] \times [-1, 1]$, the control inputs $u_{i} \in \mathcal{U}$ were chosen in the range $u_{i} \in [-1.1, 1.1]$, and the resulting states were generated according to $x_{i}{}^{\prime} \sim Q(\cdot \,|\, x_{i}, u_{i})$. 
We then presumed no knowledge of the system dynamics or the structure of the disturbance for the purpose of computing the approximately optimal control inputs using our proposed method.

Using $\mathcal{S}$, we then computed an estimate $\hat{m}$ according to \eqref{eqn: estimate closed form solution}, which can be viewed as an empirical estimate of the system dynamics. 
We used a bandwidth parameter $\sigma = 1$ for our calculations, which was determined by cross-validation.
We then chose a collection of admissible control inputs $\lbrace \tilde{u}_{\ell} \rbrace_{\ell = 1}^{L}$, $L = 100$, in the range $\tilde{u}_{\ell} \in [-1, 1]$ to compute the estimator $\hat{p}$ in \eqref{eqn: policy representation}. 
We then selected $R = 25$ \emph{evaluation points} $\lbrace x_{j} \rbrace_{j=1}^{R}$, chosen uniformly in the region $[-1, 1] \times [-1, 1]$ from which to compute the approximately optimal control inputs.

In order to demonstrate the effectiveness of the method, we computed the optimal control inputs using CVX \cite{cvx} from the evaluation points $\lbrace x_{j} \rbrace_{j=1}^{R}$ using the \emph{deterministic} dynamics. 
Once we computed the optimal inputs, we propagate the dynamics forward in time using the optimal inputs to obtain the state at the next time instant. 
We then plotted the vector field of the closed-loop dynamics under the optimal control input in Figure \ref{fig: double integrator} (blue).
We then computed the approximately optimal control inputs using \eqref{eqn: approximate optimization problem} using the sample $\mathcal{S}$ taken from the \emph{stochastic} dynamics to minimize the cost $c$ at each point $x_{j}$ over a single time step. Using the computed inputs, we then plotted the vector field of the closed-loop dynamics in Figure \ref{fig: double integrator} (orange) to compare against the optimal control inputs computed via CVX. 

\begin{figure}[!t]
    \centering
    \includegraphics[width=240pt, height=120pt]{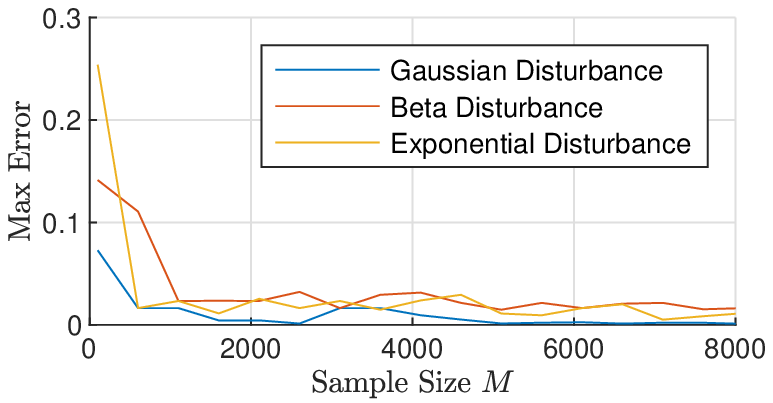}
    \includegraphics[width=240pt, height=120pt]{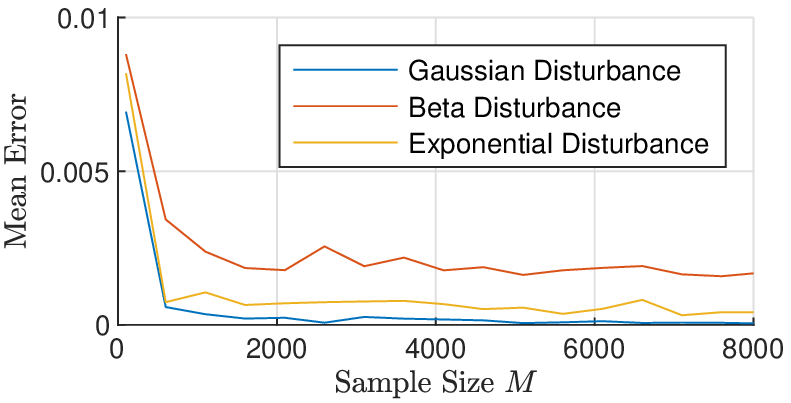}
    \caption{%
        (Top) Maximum error of the control inputs computed via our proposed method vs. the optimal control inputs computed via CVX. 
        (Bottom) Mean error of the control inputs.
        The mean error is less than $0.005$ when $M > 1000$. 
    }
    \label{fig: double integrator mean error}
\end{figure}

We can see in Figure \ref{fig: double integrator} that the control inputs selected by the algorithm are close to the optimal control inputs obtained by CVX, especially further away from the origin, where the computed control inputs coincide almost exactly with the optimal control inputs. 
Closer to the origin, we see that the algorithm deviates slightly from the optimal control inputs, which we anticipate is due largely to the randomness of the state observations.

We then generated a new sample $\mathcal{S}$ of the system in \eqref{eqn: double integrator dynamics} affected by a disturbance with a beta distribution, and then from \eqref{eqn: double integrator dynamics} affected by a disturbance with an exponential distribution and computed the optimal control inputs using our proposed method. 
The vector fields of the closed-loop dynamics for these cases are shown in Figure \ref{fig: double integrator non gaussian}. 
We can see that the algorithm computes an approximately optimal controller, despite the state observations being affected by a non-Gaussian disturbance. 

\begin{figure}[!t]
    \centering
    \includegraphics[width=240pt, height=120pt]{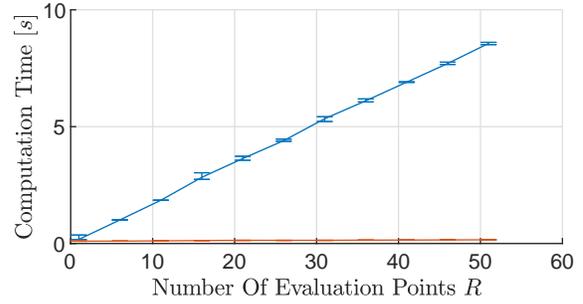}
    \caption{%
        Mean computation time of the optimal control solution, computed using CVX \cite{cvx} (blue) vs. the mean computation time of the kernel based algorithm (orange) as a function of the number of evaluation points $R$. 
        The sample size used to construct the estimate $\hat{m}$ is $M = 1600$, and the number of admissible control inputs is $P = 100$.
    }
    \label{fig: computation time function of R}
\end{figure}

\begin{figure}[!t]
    \centering
    \includegraphics[width=240pt, height=120pt]{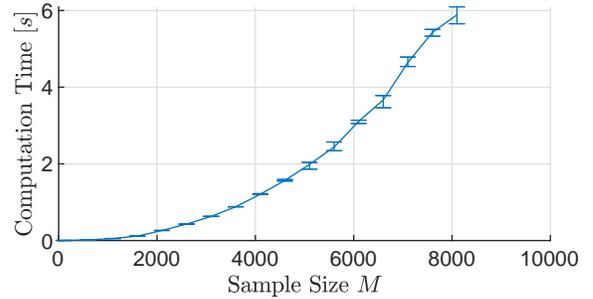}
    \caption{%
        Mean computation time of the kernel based algorithm (blue) as a function of the sample size $M$ used to construct the approximation. The number of evaluation points is $R = 10$, and the number of admissible control inputs is $P = 100$. The computation time increases exponentially as the sample size increases.
    }
    \label{fig: computation time function of M}
\end{figure}

Note that the quality of the approximation obtained from our method depends on the sample sizes $M$. 
As the number of observations in the sample increase, the approximately optimal control inputs converge to the actual optimal control inputs. 
In order to demonstrate this, we computed the maximum and mean squared Euclidean error between the control inputs computed via our method and the optimal control inputs computed via CVX for varying sample sizes $M \in [100, 8000]$ in order to characterize the performance of the algorithm. The results are shown in Figure \ref{fig: double integrator mean error}. We can see that the error of the approximately optimal control inputs decreases quickly as the sample size increases, and that the mean error is approximately less than $0.005$ when the sample size is greater than $M > 2000$. 
However, we can also see that the quality of the approximation does not improve significantly as the sample size increases, which is due to the asymptotic convergence of the estimate $\hat{m}$ to the true embedding $m$ \cite{song2009hilbert}. This presents a tradeoff between computation time and numerical accuracy, especially since the complexity scales exponentially as the sample size increases. 

The computation times for both approaches are shown in Figure \ref{fig: computation time function of R} as a function of the evaluation points up to $R = 51$.
We can see from Figure \ref{fig: computation time function of R} that the computation times for the CVX optimization method increases roughly linearly as the number of evaluation points increases, since the optimization problem needs to solve for each point independently. The computation time for our proposed method is also roughly linear in the number of evaluation points, but is dominated primarily by the computation time required for the matrix inversion, which increases exponentially with the sample size $M$, and is generally $\mathcal{O}(M^{3})$ \cite{song2009hilbert}.

This is demonstrated empirically in Figure 
\ref{fig: computation time function of M}, 
where we compute the mean computation time as a function of the sample size $M$.
We can see that as the sample size increases, the computation time increases exponentially.
However, as mentioned earlier, the computation time of the kernel based approach can also be improved using existing speedup techniques \cite{rahimi2007random, le2013fastfood}.

This also illustrates the computational advantage of our proposed method, since the quality of the approximation obtained via kernel methods has a mean error of roughly $0.005$ with a sample size of $M = 1600$, but is able to compute the optimal control inputs an order of magnitude faster than the optimal solution via CVX for multiple evaluation points.


\subsection{Nonholonomic Vehicle}

\begin{figure}[!t]
    \centering
    \includegraphics[width=240pt, height=180pt]{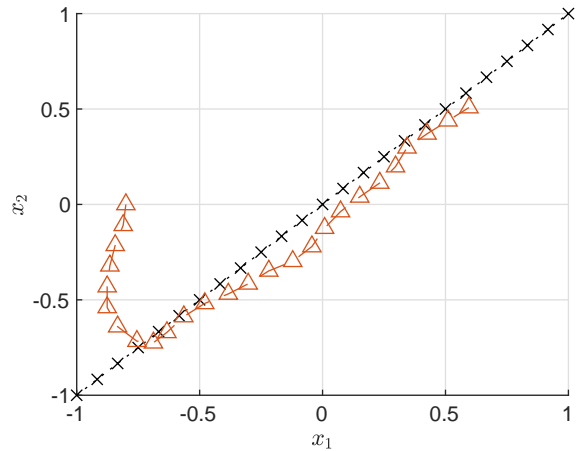}
    \caption{%
        Trajectory computed using our proposed method (orange) which tracks the target trajectory (black). The control actions are computed forward in time, with the system selecting the approximately optimal control action at each time step.
    }
    \label{fig: nonholonomic vehicle tracking forward}
\end{figure}

\begin{figure}[!t]
    \centering
    \includegraphics[width=240pt, height=180pt]{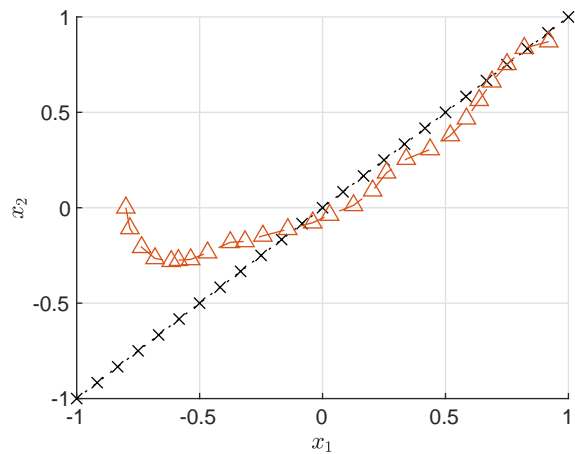}
    \caption{%
        Trajectory computed using our proposed method (orange) which tracks the target trajectory (black). The control actions are computed backward in time using dynamic programming. Note that the trajectory more closely follows the target using dynamic programming.
    }
    \label{fig: nonholonomic vehicle tracking dynamic programming}
\end{figure}

We consider the problem of target tracking for a system with nonholonomic vehicle dynamics as defined in \cite{chiang2015path}, modified such that it has a minimum forward velocity. The dynamics with sampling time $T_{s}$ are given by:
\begin{align}
    \label{eqn: nonholonomic vehicle dynamics}
    \begin{aligned}[c]
        \dot{x}_{1} &= (u_{1} + V_{\rm min}) \sin(x_{3}) + w \\ 
        \dot{x}_{2} &= (u_{1} + V_{\rm min}) \cos(x_{3}) + w \\ 
        \dot{x}_{3} &= u_{2} + w
    \end{aligned}
\end{align}
where $V_{\rm min} = 0.1$ is the minimum, constant forward velocity, $[x_{1}, x_{2}, x_{3}] \in \mathcal{X} \subseteq \mathbb{R}^{3}$ are the states, $[u_{1}, u_{2}]^{\top} \in \mathcal{U} \in \mathbb{R}^{2}$ are the control inputs, and $w \sim \mathcal{N}(0, \Sigma)$ is a random variable on the measurable space $(\mathbb{R}^{p}, \mathscr{B}(\mathbb{R}^{p}))$, where $\Sigma = 0.1 I$. 
We define a target trajectory, moving from $x_{\rm init} = [-1, -1, \pi/4]^{\top}$ to $x_{\rm final} = [1, 1, \pi/4]^{\top}$ (shown in black in Fig. \ref{fig: nonholonomic vehicle tracking forward} and Fig. \ref{fig: nonholonomic vehicle tracking dynamic programming}), and define the cost function such that the goal is to minimize the squared Euclidean distance from the system's position to the target trajectory's position at each time step.

We consider a sample $\mathcal{S} = \lbrace (x_{i}, u_{i}, x_{i}{}^{\prime}) \rbrace_{i=1}^{M}$ of size $M = 1600$ taken i.i.d. from $Q$, a representation of \eqref{eqn: nonholonomic vehicle dynamics} as a Markov control process. The states $x_{i}$ were drawn uniformly in the range $[-1.1, 1.1] \times [-1.1, 1.1] \times [-6, 6]$, the control inputs $u_{i}$ were drawn uniformly in the range $[-0.1, 1.2] \times [-10.1, 10.1]$, and then $x_{i}{}^{\prime}$ drawn from $Q(\cdot \,|\, x_{i}, u_{i})$. 

We then computed an estimate $\hat{m}$ according to \eqref{eqn: estimate closed form solution} and used a bandwidth parameter of $\sigma = 3$ for the kernel function, determined by cross-validation.
We then selected a collection $\lbrace \tilde{u}_{\ell} \rbrace_{\ell=1}^{L}$ of $L = 231$ admissible control inputs within the range $[0, 1] \times [-10, 10]$ to construct $\hat{p}$ as in \eqref{eqn: policy representation}.
We choose an initial condition $x_{0} = [-0.8, 0, \pi]^{\top}$, and evolve the system forward in time via \eqref{eqn: nonholonomic vehicle dynamics} over a time horizon $N = 20$, computing an approximately optimal control action at each time step using our proposed method.
The resulting trajectory is plotted in Figure \ref{fig: nonholonomic vehicle tracking forward} (orange), and the computation time over the time horizon was approximately $0.538$ seconds. 
As expected, we can see that the control actions selected from our proposed method drive the system to closely follow the target trajectory. 


We then computed the optimal controller via dynamic programming in order to compare against the forward in time approach. Unlike the previous approach, in which the control actions are selected in a greedy fashion, the dynamic programming approach computes the optimal control actions backward in time by iteratively optimizing a sequence of value functions \eqref{eqn: value function recursion}, and then selecting the control actions at each time step which have the highest \emph{value}. 
We use the same tracking trajectory as before, as well as the same initial condition in order to compare the performance of the two approaches.
The resulting trajectory is shown in Figure \ref{fig: nonholonomic vehicle tracking dynamic programming} (orange). 
The computation time for the dynamic programming solution was approximately $6.448$ seconds.

As expected, we can see that the trajectories obtained from the two approaches both follow the target trajectory, but the dynamic programming solution follows the trajectory better over the entire time horizon. This is because the value functions take into account the future actions of the system in order to minimize the total cost. 
This shows that our algorithm is able to select the approximately optimal control actions at each time step for a nonlinear system either forward in time or backward in time via dynamic programming, using only sample information taken from observations of the system evolution.


\section{Conclusions \& Future Work}
\label{section: concluding remarks}

In this paper, we have presented a novel method for computing the optimal policy for discrete-time dynamic programming problems using observations taken from a stochastic system under an arbitrary disturbance. Our method is model-free and largely agnostic to the cost function used. 
We have demonstrated our proposed method on a discrete time stochastic double integrator system and on a nonholonomic vehicle target tracking problem.
We plan to explore further theoretical extensions of this method to other classes of stochastic control problems, and to constrained optimal control problems.


\bibliographystyle{IEEEtran}
\bibliography{IEEEabrv, shortIEEE, bib-sorted}

\end{document}